\documentclass[twoside]{article}
\usepackage[english]{babel}
\usepackage{amssymb,amsmath,amsthm}
\begin{document}


\begin{center}
{\Large \textbf{Complete lists of low dimensional complex associative algebras}}\\
\end{center}

\begin{center}
\textbf{I.S. Rakhimov$^{1}$, I.M. Rikhsiboev$^{2}$ and W.Basri$^{3}$}\\
\end{center}
\begin{center}
$^{1,3}$Department of Mathematics, Faculty of Science, \\
$^{1,2,3}$Laboratory of Theoretical studies,
Institute for Mathematical Research (INSPEM),\\
University Putra Malaysia\\

$^{1}$\emph{risamiddin@gmail.com, $^{2}$ikromr@gmail.com, $^{3}$witri@science.upm.edu.my}

\end{center}
\thispagestyle{empty}
\begin{abstract}
In this paper we present a complete classification (isomorphism classes with some isomorphism invariants) of complex associative algebras up to dimension five (including both cases: unitary and non-unitary). In some symbolic computations we used Maple software.
\end{abstract}


\textbf{}

\section{Introduction}

\textbf{}

 In algebra, there are three strongly related classical algebras: associative, Lie and Jordan algebras. The objects of our attention in this paper are associative algebras, exactly finite dimensional associative algebras over complex numbers. The major theorems on associative algebras include the most splendid results of the great heroes of algebra: Wedderburn, Artin, Noether, Hasse, Brauer, Albert, Jacobson, and many others.

 It is known that the universal enveloping algebra of a Lie algebra has the structure of an associative algebra. Introduced by Loday \cite{LL}, the notion of Leibniz algebra is a generalization of the Lie algebra, where the skew-symmetry in the bracket is dropped. Loday \cite{LF} also showed that the relationship between Lie and associative algebras can be translated into an analogous relationship between Leibniz and the so-called dialgebras (or diassociative algebras), which are a generalization of associative algebra possessing two operations $\dashv$ and $\vdash.$ In particular, it was shown that any dialgebra becomes a Leibniz algebra under the bracket $[x,y]=x\dashv y - y\vdash x.$

The classification of low dimensional complex Lie and Leibniz algebras can be found in \cite{Jc}, \cite{LL}, \cite{AO}.

Our motivation to classify associative algebras comes out from the intention to classify the diassociative algebras. A diassociative algebra is a generalization of associative algebra, it can be obtained as a combination of two associative algebras. In order to make up this combination we need the list of all associative algebras (both unitary and non-unitary). The latest lists of all unitary associative algebras in dimension two, three, four, and five are available in \cite{PB}, \cite{Az}, \cite{Gb} and \cite{Mz}, respectively.

In 1993 Loday \cite{LL} introduced a non-antisymmetric version of Lie algebras, whose bracket satisfies the Leibniz identity and therefore this generalization has been called Leibniz algebras. The Leibniz identity, combined with antisymmetricity, is a variation of the Jacobi identity. Hence a Lie algebra is antisymmetric Leibniz algebra. Loday also introduced an ``associative'' version of Leibniz algebras, called diassociative algebras, equipped with two binary operations, which satisfy certain identities. These identities are all variations of the associative law. An associative algebra is a diassociative algebra when the two operations coincide. The main motivation of Loday to introduce this class of algebras was the search of an ``obstruction'' to the periodicity of algebraic K-theory. Besides this purely algebraic motivation, some relationships with classical geometry, non-commutative geometry and physics have been recently discovered. We are interested in classification of diassociative algebras. One of the approaches to describe a finite dimensional diassociative law is to consider it as a combination of two associative algebras. Therefore, we need a complete list of associative algebras. Obviously, the classification problem of algebras (even associative algebras), in general, is nearly unreal problem.  All existence classifications of associative algebras concern unitary ones. In order to use the classification of associative algebras to classify the diassociative algebras we need both lists of unitary and non unitary associative algebras. In the present, paper we give complete lists of low-dimensional complex associative algebras.

\section{Lists of low-dimensional associative algebras with some isomorphism invariants}
Now and what it follows, all algebras are assumed to be over the field of complex numbers $\mathbb{C}$. \\
In the next section we give lists of all complex associative algebras in dimensions 2-- 4. Some remarks on the tables. In the tables $As_{p}^{q} $ stands for $q^{th}$ algebra in dimension $p$. In the second column only nonzero products are given. The third column describes the automorphism groups of the algebras.  For all algebras we test to be nilpotent or not, in the last case we indicate nilpotent $N$ and semi-simple $S$ parts of $As_{p}^{q} $. By $C(As_{p}^{q}),$ $L(As_{p}^{q})$, and $R(As_{p}^{q})$ the maximal commutative subalgebra, the left annihilator, and
the right annihilator of $As_{p}^{q}$, respectively are denoted, and their dimensions are given in the last three columns of the tables.

A proof, that any associative algebra of dimensions 2-- 4 is included in the lists, is available from the authors. Because of its length, it is omitted from the paper. We tabulate only indecomposable algebras.

\subsection{Two-dimensional associative algebras}
\begin{tabular}{|c|p{1.4in}|p{1in}|p{1.2in}|p{0.4in}|p{0.4in}|p{0.4in}|} \hline
\footnotesize\textit{} & \textit{Table of multiplication} & \textit{ Automorphisms} & \textit{Type of algebra}&\textit{dim} & \textit{dim} & \textit{dim} \\

&& & &\textit{C$(As_{p}^{q})$}&\textit{L$(As_{p}^{q})$}& \textit{R$(As_{p}^{q})$} \\
\hline
$As_{2}^{1}:$&${e_{1}e_{1}=e_{2}}$ & $\, \left(\begin{array}{l} {a\, \, \, \, \, 0} \\ {b\, \, \, \, a^{2} }
\end{array}\right)$& \textit{commutative,} \newline \textit{nilpotent}& {2} & {1} & {1} \\\hline
$As_{2}^{2}:$&$e_{1}e_{1} =e_{1},$ $e_{1}e_{2}=e_{2}$ & $\, \left(\begin{array}{l} {1\, \, \, \, 0} \\ {a\, \, \, \, b} \end{array}\right)$ & $A=N \dotplus S,$ \newline $N=<e_2>$, \newline $S=<e_1>$& 1 & 1 & 0 \\\hline
$As_{2}^{3}:$&$e_{1}e_{1}=e_{1},$ $e_{2}e_{1}=e_{2}$ & $\, \left(\begin{array}{l} {1\, \, \, \, 0} \\ {a\, \, \, \, b} \end{array}\right)$ & $A=N \dotplus S,$  \newline $N=<e_2>$, \newline $S=<e_1>$ & 1 & 0 & 1 \\\hline
$As_{2}^{4}:$&$e_{1}e_{1}=e_{1},$ $e_{1}e_{2}=e_{2},$ \newline $e_{2} e_{1}=e_{2} $& $\left(\begin{array}{l} {1\, \, \, \, 0} \\ {0\, \, \, \, a} \end{array}\right)$ & \textit{commutative,} \newline \textit{unitary}, \newline $A=N \dotplus S,$ \newline $N=<e_2>$,\newline $S=<e_1>$ & 2 & 0 & 0 \\\hline
\end{tabular}

\subsection{Three-dimensional associative algebras}

\begin{tabular}{|c|p{1.3in}|p{1.8in}|p{1in}|p{0.35in}|p{0.35in}|p{0.35in}|} \hline
&\textit{Table of} \newline \textit {multiplication}& \textit{ Automorphisms} & \textit{Type of algebra} & \textit {dim C$(As_{p}^{q})$} & \textit{dim L$(As_{p}^{q})$} & \textit{dim R$(As_{p}^{q})$} \\

\hline
$As_{3}^{1}:$ & $e_{1}e_{3}=e_{2},$ $e_{3}e_{1}=e_{2}$ & $\left(\begin{array}{l} {a\, \, \, \, \, 0\, \, \, \, \, 0} \\ {b\, \, \, \, ac\, \, \, d} \\ {0\, \, \, \, \, 0\, \, \, \, \, c} \end{array}\right)$, $\left(\begin{array}{l} {0\, \, \, \, \, 0\, \, \, \, \, a} \\ {b\, \, \, \, ac\, \, \, d} \\ {c\, \, \, \, \, 0\, \, \, \, \, 0} \end{array}\right)$ & \textit{commutative,} \newline \textit{nilpotent}& 3 & 1 & 1 \\ \hline
$As_{3}^{2}:$ & $e_{1}e_{3}=e_{2},$ \newline $ e_{3}e_{1}=\alpha e_{2},$ \newline $\alpha$ $\in \mathbb{C} \backslash\ \{1\}$ & $\, \left(\begin{array}{l} {a\, \, \, \, \, 0\, \, \, \, \, 0} \\ {b\, \, \, \, ac\, \, \, d} \\ {0\, \, \, \, \, 0\, \, \, \, \, c} \end{array}\right)$,$\, \left(\begin{array}{l} {0\, \, \, \, \, 0\, \, \, \, \, a} \\ {b\, \, \, \, ac\, \, \, d} \\ {c\, \, \, \, \, 0\, \, \, \, \, 0} \end{array}\right)$,\newline  $\left(\begin{array}{l} {a\, \, \, \, \, \, \, \, \, \, 0\, \, \, \, \, \, \, \, \, \, b} \\ {c\, \, \, \, af-be\, \, \, d} \\ {e\, \, \, \, \, \, \, \, \, \, 0\, \, \, \, \, \, \, \, \, \, f} \end{array}\right)$ & \textit{nilpotent} & 2 & 1 & 1 \\ \hline
$As_{3}^{3}:$ & $e_{1}e_{1}=e_{2},$ $e_{1}e_{2}=e_{3}$, \newline $e_{2}e_{1}=e_{3}$ & $ \left(\begin{array}{l} {a\, \, \, \, \, \, \, 0\, \, \, \, \, \, 0} \\ {b\, \, \, \, \, \, a^{2} \, \, \, \, \, 0} \\ {c\, \, \, 2ab\, \, \, a^{3} } \end{array}\right)$ & \textit{commutative,} \newline \textit{nilpotent} & 3 & 1 & 1 \\
\hline
\end{tabular}

\begin{tabular}{|c|p{1.3in}|p{1.5in}|p{1in}|p{0.35in}|p{0.35in}|p{0.35in}|} \hline
$As_{3}^{4}:$ & $e_{1}e_{3}=e_{2},$ $e_{2}e_{3}=e_{2},$ \newline $e_{3}e_{3}=e_{3}$ & $\, \left(\begin{array}{l} {a\, \, \, \, \, \, \, 0\, \, \, \, \, \, \, \, 0} \\ {b\, \, \, \, a+b\, \, \, c} \\ {0\, \, \, \, \, \, \, \, 0\, \, \, \, \, \, \, \, 1} \end{array}\right)$ & $A=N \dotplus S, $ \newline $N=<e_1, e_2>,$ \newline $S=<e_3>$& 2 & 0 & 2 \\
\hline
$As_{3}^{5} :$ & $e_{2}e_{3}=e_{2},$ $e_{3}e_{1}=e_{1},$ \newline $e_{3}e_{3}=e_{3}$ & $\, \left(\begin{array}{l} {a\, \, \, \, 0\, \, \, b} \\ {0\, \, \, \, c\, \, \, \, d} \\ {0\, \, \, \, 0\, \, \, \, 1} \end{array}\right)$ & $A=N \dotplus S, $\newline $N=<e_1, e_2>,$ \newline $S=<e_3>$ & 2 & 1 & 1 \\ \hline
$As_{3}^{6}:$ & $e_{3}e_{1}=e_{2},$  $e_{3}e_{2}=e_{2},$ \newline $e_{3}e_{3}=e_{3}$ & $\, \left(\begin{array}{l} {a\, \, \, \, \, \, \, 0\, \, \, \, \, \, \, \, 0} \\ {b\, \, \, \, a+b\, \, \, c} \\ {0\, \, \, \, \, \, \, \, 0\, \, \, \, \, \, \, \, 1} \end{array}\right)$ & $A=N \dotplus S, $ \newline $N=<e_1, e_2>,$ \newline $S=<e_3>$& 2 & 2 & 0 \\
\hline
$As_{3}^{7} :$ & $e_{1}e_{2}=e_{1},$ $e_{2}e_{2}=e_{2},$ \newline $e_{3}e_{1}=e_{1},$ $e_{3}e_{3}=e_{3} $ & $\, \left(\begin{array}{l} {a\, \, \, \, b\, \, \, -b} \\ {0\, \, \, \, \, 1\, \, \, \, \, \, \, 0} \\ {0\, \, \, \, 0\, \, \, \, \, \, \, \, 1} \end{array}\right)$ & \textit{unitary},\newline $A=N \dotplus S, $ \newline $N=<e_1>,$ \newline $S=<e_2, e_3>$ & 2 & 0 & 0 \\ \hline
$As_{3}^{8}:$ & $e_{1}e_{3}=e_{1},$  $e_{2}e_{3}=e_{2},$ \newline $e_{3}e_{1}=e_{1},$  $e_{3}e_{3}=e_{3}$  & $\, \left(\begin{array}{l} {a\, \, \, \, 0\, \, \, \, 0} \\ {0\, \, \, \, b\, \, \, \, c} \\ {0\, \, \, \, 0\, \, \, \, 1} \end{array}\right)$ & $A=N \dotplus S, $ \newline $N=<e_1, e_2>,$ \newline  $S=<e_3>$ & 2 & 0 & 1 \\ \hline
$\, As_{3}^{9}:$ \, & $e_{2}e_{3}=e_{2},$ $e_{3}e_{1}=e_{1},$ \newline $e_{3}e_{2}=e_{2},$ $e_{3}e_{3}=e_{3}$ & $\, \left(\begin{array}{l} {a\, \, \, \, 0\, \, \, \, b} \\ {0\, \, \, \, c\, \, \, \, 0} \\ {0\, \, \, \, 0\, \, \, \, 1} \end{array}\right)$ & $A=N \dotplus S, $ \newline $N=<e_1, e_2>,$ \newline $S=<e_3>$ & 2 & 1 & 0 \\ \hline
$\, As_{3}^{10}:\, $ & $e_{1}e_{3}=e_{1},$ $e_{2}e_{3}=e_{2},$ \newline $e_{3}e_{1}=e_{1},$ $e_{3}e_{2}=e_{2},$ \newline $e_{3}e_{3}=e_{3}$ & $\, \left(\begin{array}{l} {a\, \, \, \, b\, \, \, \, 0} \\ {c\, \, \, \, d\, \, \, \, 0} \\ {0\, \, \, \, 0\, \, \, \, 1} \end{array}\right)$ & \textit{commutative, unitary,} \newline $A=N \dotplus S,$ \newline $N=<e_1, e_2>,$ \newline $S=<e_3>$ & 3 & 0 & 0 \\ \hline
$\, As_{3}^{11}:$ \, & $e_{1}e_{3}=e_{2},$ $e_{2} e_{3}=e_{2},$ \newline $e_{3} e_{1}=e_{2},$ $e_{3} e_{2}=e_{2},$ \newline $e_{3} e_{3}=e_{3}$ & $\, \left(\begin{array}{l} {a\, \, \, \, \, \, \, 0\, \, \, \, \, \, \, \, 0} \\ {b\, \, \, \, a+b\, \, \, 0} \\ {0\, \, \, \, \, \, \, \, 0\, \, \, \, \, \, \, \, 1} \end{array}\right)$ & \textit{commutative,}\newline $A=N \dotplus S, $ \newline $N=<e_1, e_2>,$ \newline $S=<e_3>$& 3 & 0 & 0 \\ \hline
$\, As_{3}^{12}: \,$ & $e_{1}e_{1}=e_{2},$ $e_{1}e_{3}=e_{1},$ \newline $e_{2}e_{3}=e_{2},$ $e_{3}e_{1}=e_{1},$ \newline $e_{3}e_{2}=e_{2},$  $e_{3}e_{3}=e_{3}$ & $\, \left(\begin{array}{l} {a\, \, \, \, \, 0\, \, \, \, \, 0} \\ {b\, \, \, \, \, a^{2} \, \, \, 0} \\ {0\, \, \, \, \, \, 0\, \, \, \, \, \, 1} \end{array}\right)$ & \textit{commutative, unitary,} \newline $A=N \dotplus S,$ \newline $N=<e_1, e_2>,$ \newline $S=<e_3>$ & 3 & 0 & 0 \\ \hline
\end{tabular}

\subsection{Four-dimensional associative algebras}

\begin{footnotesize} 
\begin{tabular}{|c|p{1.1in}|p{1.9in}|p{1in}|p{0.3in}|p{0.3in}|p{0.3in}|} \hline
& \textit{Table of} \newline \textit{multiplication}& \textit{Automorphisms} & \textit{Type of algebra} & \textit{dim C$(As_{p}^{q})$} & \textit{dim L$(As_{p}^{q})$} & \textit{dim R$(As_{p}^{q})$} \\
\hline
$As_{4}^{1}:$ &$e_{1}e_{2}=e_{3},$ $e_{2}e_{1}=e_{4}$ & $\left(\begin{array}{l} {a\, \, \, \, 0\, \, \, \, 0\, \, \, \, \, 0} \\ {0\, \, \, \, a\, \, \, \, 0\, \, \, \, \, 0} \\ {b\, \, \, \, c\, \, \, \, a^{2} \, \, \, 0} \\ {d\, \, \, \, e\, \, \, \, 0\, \, \, \, a^{2} } \end{array}\right)$, $\left(\begin{array}{l} {0\, \, \, \, a\, \, \, \, 0\, \, \, \, \, 0} \\ {b\, \, \, \, 0\, \, \, \, 0\, \, \, \, \, 0} \\ {c\, \, \, \, d\, \, \, \, 0\, \, \, ab} \\ {e\, \, \, \, f\, \, \, \, ab\, \, \, \, 0} \end{array}\right)$ & \textit{nilpotent} & 3 & 2 & 2 \\
\hline
$As_{4}^{2}:$ &$e_{1}e_{2}=e_{4},$ $e_{3}e_{1}=e_{4}$ & $\left(\begin{array}{l} {\, \, a\, \, \, \, \, 0\, \, \, \, \, 0\, \, \, \, \, 0} \\ {\, \, c\, \, \, \, \, b\, \, \, \, \, 0\, \, \, \, \, 0} \\ {-c\, \, \, 0\, \, \, \, b\, \, \, \, \, 0} \\ {\, \, d\, \, \, \, e\, \, \, \, f\, \, \, \, ab} \end{array}\right)$& \textit{nilpotent} & 3 & 2 & 2 \\ \hline
$As_{4}^{3}:$ & $e_1e_2=e_3$, $e_2e_1=e_4$, \newline$e_2e_2=-e_3$ & $\left(\begin{array}{l} {a\, \, \, \, 0\, \, \, \, 0\, \, \, \, \, 0} \\ {0\, \, \, \, a\, \, \, \, 0\, \, \, \, \, 0} \\ {b\, \, \, \, c\, \, \, \, a^{2} \, \, \, 0} \\ {d\, \, \, \, e\, \, \, \, 0\, \, \, \, a^{2} } \end{array}\right)$\textit{} & \textit{nilpotent} & 3 & 2 & 2 \\ \hline
$As_{4}^{4}:$& $e_1e_2=e_3,$ $e_2e_2=e_4$, \newline $e_2e_1=-e_3$ & $\left(\begin{array}{l} {a\, \, \, \, c\, \, \, \, 0\, \, \, \, \, 0} \\ {0\, \, \, \, b\, \, \, \, 0\, \, \, \, \, 0} \\ {d\, \, \, \, e\, \, \, \, ab\, \, \, 0} \\ {f\, \, \, \, g\, \, \, \, 0\, \, \, \, b^{2} } \end{array}\right)$\textit{} & \textit{nilpotent} & 3 & 2 & 2 \\
\hline
$As_{4}^{5}:$& $e_1e_2=e_4,$ $e_3e_3=e_4$ \newline $e_2e_1=-e_4,$ & $\left(\begin{array}{l} {a\, \, -\frac{b^{2} }{c} \, \, \, \, 0\, \, \, \, \, 0} \\ {c\, \, \, \, \, \, \, 0\, \, \, \, \, \, \, 0\, \, \, \, \, \, \, 0} \\ {0\, \, \, \, \, \, \, 0\, \, \, \, \, \, \, b\, \, \, \, \, \, \, 0} \\ {d\, \, \, \, \, \, \, e\, \, \, \, \, \, \, f\, \, \, \, \, b^{2} } \end{array}\right)$, \textit{} & \textit{nilpotent} & 3 & 1 & 1 \\
&&$\left(\begin{array}{l} {\frac{d^{2}+ab}{c}\, \, \, a \, \, \, \, 0\, \, \, \, \, 0} \\ {\, \, \, \, b\, \, \, \, \, \, \, \, c\, \, \, \, \, \, \, 0\, \, \, \, \, \, \, 0} \\ {\, \, \, \, 0\, \, \, \, \, \, \, \, 0\, \, \, \, \, \, \, d\, \, \, \, \, \, \, 0} \\ {\, \, \, \, e\, \, \, \, \, \, \, \, f\, \, \, \, \, \, \, h\, \, \, \, \, d^{2} } \end{array}\right)$&&&& \\ \hline
$As_{4}^{6}(\alpha ):$ &$e_1e_2=e_4$, $e_2e_2=e_3,$ \newline $e_2e_1=\frac{1+\alpha}{1-\alpha}e_4$, $\alpha\neq 1$ & $\left(\begin{array}{l} {a\, \, \, \, c\, \, \, \, \, \, \, \, \, \, \, \, 0\, \, \, \, \, \, \, \, \, \, \, \, \, \, \, \, \, \, 0} \\ {0\, \, \, \, b\, \, \, \, \, \, \, \, \, \, \, \, 0\, \, \, \, \, \, \, \, \, \, \, \, \, \, \, \, \, \, 0} \\ {d\, \, \, \, e\, \, \, \, \, \, \, \, \, \, \, \, b^{2} \, \, \, \, \, \, \, \, \, \, \, \, \, \, \, 0} \\ {f\, \, \, \, g\, \, \, \, bc(1+\alpha )\, \, \, \, ab}\end{array}\right)$ & \textit{nilpotent} & 3 & 2 & 2 \\\hline
$As_{4}^{7}:$ & $e_1e_1=e_1,$ $e_1e_4=e_4,$ \newline $e_2e_1=e_2,$ $e_2e_4=e_3 $& $\left(\begin{array}{l} {1\, \, \, \, \, 0\, \, \, \, \, 0\, \, \, \, \, 0} \\ {a\, \, \, \, \, b\, \, \, \, \, 0\, \, \, \, \, 0} \\ {ac\, \, bc\, \, bd\, \, ad} \\ {c\, \, \, \, \, 0\, \, \, \, \, 0\, \, \, \, \, d} \end{array}\right) $ & $A=N \dotplus S, $ \newline $N=<e_2, e_3, e_4 >$, \newline $S=<e_1>$ & 2 & 2 & 2 \\ \hline
$As_{4}^{8}:$ & $e_1e_1=e_1,$ $e_3e_1=e_3,$ \newline $e_4e_3=e_2,$ $e_2e_1=e_2$ & $\left(\begin{array}{l} {\, \, 1\, \, \, \, \, \, \, \, 0\, \, \, \, \, 0\, \, \, \, 0} \\ {\, \, a\, \, \, \, \, \, bc\, \, \, \, d\, \, \, \, e} \\ {-\frac{e}{c} \, \, \, 0\, \, \, \, \, b\, \, \, \, \, 0} \\ {\, \, 0\, \, \, \, \, \, 0\, \, \, \, \, 0\, \, \, \, c} \end{array}\right)$ & $A=N \dotplus S, $ \newline $N=<e_2, e_3, e_4 >$, \newline $S=<e_1>$ & 2 & 0 & 2 \\ \hline
$As_{4}^{9}:$ & $e_1e_1=e_1,$ $e_1e_3=e_3,$ \newline $e_1e_2=e_2,$ $e_3e_4=e_2$ & $ \left(\begin{array}{l} {1\, \, \, \, \, 0\, \, \, \, \, 0\, \, \, \, \, \, 0} \\ {a\, \, \, \, bc\, \, \, d\, \, \, -ce} \\ {e\, \, \, \, \, 0\, \, \, \, \, b\, \, \, \, \, \, 0} \\ {0\, \, \, \, \, 0\, \, \, \, \, 0\, \, \, \, \, \, c} \end{array}\right)$& $A=N \dotplus S, $ \newline $N=<e_2, e_3, e_4 >$, \newline $S=<e_1>$ & 2 & 2 & 0 \\ \hline
$As_{4}^{10}:$ & $e_1e_1=e_3,$ $e_1e_3=e_4,$ \newline $e_2e_2=-e_4,$ \newline $e_3e_1=e_4$ & $\left(\begin{array}{l} {\sqrt[{3}]{a^{2} } \, \, 0\, \, \, \, \, \, \, \, \, \, \, \, \, 0\, \, \, \, \, \, \, \, \, \, \, \, \, \, \, \, \, 0} \\ {\, b\, \, \, \, \, \, \, \, a\, \, \, \, \, \, \, \, \, \, \, \, \, \, 0\, \, \, \, \, \, \, \, \, \, \, \, \, \, \, \, \, 0} \\ {\, c\, \, \, \, \, \, b\sqrt[{3}]{a} \, \, \, \, \, \, \, \, \, \sqrt[{3}]{a^{2} } \, \, \, \, \, \, \, \, \, 0} \\ {d\, \, \, \, \, \, e\, \, \, \, \, \, 2c\sqrt[{3}]{a^{2} } -b^{2} \, \, \, a^{2} } \end{array}\right)$ & \textit{commutative,} \newline \textit{nilpotent} & 4 & 1 & 1 \\ \hline
$As_{4}^{11}:$ & $e_1e_1=e_4,$ $e_2e_1=e_3,$ \newline $e_1e_4=-e_3,$\newline $e_4e_1=-e_3$ & $ \left(\begin{array}{l} {a\, \, \, \, \, 0\, \, \, \, \, 0\, \, \, \, \, \, \, \, \, \, \, 0} \\ {b\, \, \, \, a^{2} \, \, \, 0\, \, \, \, \, \, \, \, \, \, \, 0} \\ {c\, \, \, \, \, d\, \, \, \, \, a^{3} \, \, a(b-2e)} \\ {e\, \, \, \, \, 0\, \, \, \, \, 0\, \, \, \, \, \, \, \, \, \, \, a^{2}} \end{array}\right)$  & \textit{nilpotent} & 3 & 1 & 2 \\ \hline
$As_{4}^{12}:$ & $e_1e_1=e_1,$ $e_1e_2=e_2,$ \newline $e_3e_1=e_3,$  $e_4e_1=e_4$ & $\left(\begin{array}{l} {1\, \, \, \, \, 0\, \, \, \, 0\, \, \, \, \, 0} \\ {a\, \, \, \, b\, \, \, \, 0\, \, \, \, \, 0} \\ {c\, \, \, \, \, 0\, \, \, \, d\, \, \, \, e} \\ {f\, \, \, \, 0\, \, \, \, g\, \, \, \, h} \end{array}\right)$ & $A=N \dotplus S, $ \newline $N=<e_2, e_3, e_4 >$, \newline $S=<e_1>$& 3 & 1 & 2 \\\hline
$As_{4}^{13}:$ & $e_2e_2=e_2,$  $e_2e_3=e_3,$\newline $e_2e_4=e_4,$  $e_1e_2=e_1$ & $\left(\begin{array}{l} {a\, \, \, \, \, b\, \, \, \, 0\, \, \, \, \, 0} \\ {0\, \, \, \, \, 1\, \, \, \, 0\, \, \, \, \, 0} \\ {0\, \, \, \, \, c\, \, \, \, d\, \, \, \, e} \\ {0\, \, \, \, f\, \, \, \, g\, \, \, \, h} \end{array}\right)$\textit{} & $A=N \dotplus S, $ \newline $N=<e_1, e_3, e_4 >$, \newline $S=<e_2>$ & 3& 2 & 1 \\\hline

\end{tabular}

\begin{tabular}{|c|p{1.2in}|p{1.8in}|p{1in}|p{0.25in}|p{0.25in}|p{0.25in}|} \hline
$As_{4}^{14}:$ & $e_1e_1=e_1,$  $e_3e_1=e_3,$ \newline $e_4e_1=e_4,$  $e_2e_1=e_2$ & $\left(\begin{array}{l} {1\, \, \, \, \, 0\, \, \, \, 0\, \, \, \, \, 0} \\ {a\, \, \, \, \, b\, \, \, \, c\, \, \, \, \, d} \\ {e\, \, \, \, f\, \, \, g\, \, \, \, h} \\ {i\, \, \, \, \, j\, \, \, \, k\, \, \, \, l} \end{array}\right)$\textit{} & $A=N \dotplus S, $ \newline $N=<e_2, e_3, e_4 >,$ \newline $S=<e_1>$ & 3 & 0 & 3 \\ \hline
$As_{4}^{15}:$ & $e_2e_2=e_2,$   $e_2e_3=e_3,$ \newline $e_2e_4=e_4,$  $e_2e_1=e_1$& $ \left(\begin{array}{l} {a\, \, \, \, \, b\, \, \, \, c\, \, \, \, \, d} \\ {0\, \, \, \, \, 1\, \, \, \, 0\, \, \, \, \, 0} \\ {e\, \, \, \, f\, \, \, g\, \, \, \, h} \\ {i\, \, \, \, \, j\, \, \, \, k\, \, \, \, l} \end{array}\right)$\textit{} & $A=N \dotplus S, $ \newline $N=<e_1, e_3, e_4 >$, \newline $S=<e_2>$ & 3& 3 & 0 \\ \hline
$As_{4}^{16}(\alpha ):$ & $e_1e_1=e_4,$  $e_1e_2=e_4,$ \newline $e_2e_1=\alpha e_4,$ \newline $e_3e_3=e_4$ & $\left(\begin{array}{l} {a\, \, \, 0\, \, \, 0\, \, \, \, 0} \\ {0\, \, \, a\, \, \, 0\, \, \, \, 0} \\ {0\, \, \, 0\, \, \, a\, \, \, 0} \\ {b\, \, \, c\, \, \, d\, \, \, a^{2} } \end{array}\right) $, $\left(\begin{array}{l} {a\, \, \, \, 0\, \, \, \, 0\, \, \, \, 0} \\ {b\, \, \, \frac{c^{2}}{a}\, \, \, 0\, \, \, \, 0} \\ {0\, \, \, \, 0\, \, \, c\, \, \, \, 0} \\ {d\, \, \, \, e\, \, \, f\, \, \, c^{2}} \end{array}\right) $& \textit{nilpotent} & 3 & 1 & 1 \\ \hline
$As_{4}^{17}:$& $e_1e_1=e_4,$  $e_1e_2=e_3,$ \newline $e_2e_1=-e_3,$ \newline $e_2e_2=-2e_3+e_4$ & $ \left(\begin{array}{l} {a\, \, \, \, 0\, \, \, \, 0\, \, \, \, \, 0} \\ {0\, \, \, \, a\, \, \, \, 0\, \, \, \, \, 0} \\ {b\, \, \, \, c\, \, \, \, a^{2} \, \, \, 0} \\ {d\, \, \, \, e\, \, \, \, 0\, \, \, \, a^{2} } \end{array}\right) $, $\left(\begin{array}{l} {0\, \, \, \, a\, \, \, \, \, \, 0\, \, \, \, \, \, \, \, \, \,  0} \\ {a\, \, \, \, 0\, \, \, \, \, \, 0\, \, \, \, \, \, \, \, \, \, 0} \\ {b\, \, \, \, c\, \, -a^{2}\, -2a^{2}} \\ {d\, \, \, \, e\, \, \, \, \, \,  0\, \, \, \, \, \, \, \, \, \, a^{2} } \end{array}\right) $ & \textit{nilpotent} & 3 & 2 & 2 \\ \hline
$As_{4}^{18}(\alpha):$& $e_1e_1=e_4,$  $e_1e_2=e_3,$ \newline $e_2e_1=-\alpha e_4,$ \newline $e_2e_2=-e_3$ & $ \left(\begin{array}{l} {a\, \, \, \, 0\, \, \, \, 0\, \, \, \, \, 0} \\ {0\, \, \, \, a\, \, \, \, 0\, \, \, \, \, 0} \\ {b\, \, \, \, c\, \, \, \, a^{2} \, \, \, 0} \\ {d\, \, \, \, e\, \, \, \, 0\, \, \, \, a^{2} } \end{array}\right)$, $\left(\begin{array}{l} {0\, \, \, \, a\, \, \, \, \, \, 0\, \, \, \, \, \, \, \, \, \,  0} \\ {a\, \, \, \, 0\, \, \, \, \, \, 0\, \, \, \, \, \, \, \, \, \, 0} \\ {b\, \, \, \, c\, \, -a^{2}\, -2a^{2}} \\ {d\, \, \, \, e\, \, \, \, \, \,  0\, \, \, \, \, \, \, \, \, \, a^{2} } \end{array}\right) $ & \textit{nilpotent} & 3 & 2 & 2 \\ \hline
$As_{4}^{19} :$& $e_1e_1=e_1,$ $e_2e_2=e_2,$  \newline $e_2e_3=e_3,$  $e_3e_1=e_3,$ \newline $e_4e_1=e_4$ & $\, \left(\begin{array}{l} {1\, \, \, \, \, \, 0\, \, \, \, \, 0\, \, \, \, \, 0} \\ {0\, \, \, \, \, \, 1\, \, \, \, \, 0\, \, \, \, \, 0} \\ {a\, \, \, -a\, \, \, \, b\, \, \, \, 0} \\ {c\, \, \, \, \, \, 0\, \, \, \, \, 0\, \, \, \, d} \end{array}\right)\, \, \, $& $A=N \dotplus S, $ \newline $N=<e_3, e_4>,$ \newline $S=<e_1, e_2>$ & 2 & 0 & 1 \\ \hline
$As_{4}^{20} :$& $e_1e_1=e_1,$  $e_1e_3=e_3,$ \newline $e_1e_4=e_4,$   $e_2e_2=e_2,$ \newline $e_3e_2=e_3$ & $ \left(\begin{array}{l} {1\, \, \, \, \, \, 0\, \, \, \, \, 0\, \, \, \, \, 0} \\ {0\, \, \, \, \, \, 1\, \, \, \, \, 0\, \, \, \, \, 0} \\ {a\, \, \, -a\, \, \, \, b\, \, \, \, 0} \\ {c\, \, \, \, \, \, 0\, \, \, \, \, 0\, \, \, \, d} \end{array}\right)\, \, \, $ & $A=N \dotplus S, $ \newline $N=<e_3, e_4>,$ \newline $S=<e_1, e_2>$ & 2 & 1 & 0 \\ \hline
$As_{4}^{21}:$ & $e_1e_1=e_1,$   $e_2e_2=e_2,$ \newline $e_2e_4=e_4,$   $e_4e_1=e_4,$ \newline $e_3e_2=e_3$ & $ \left(\begin{array}{l} {1\, \, \, \, \, \, 0\, \, \, \, \, 0\, \, \, \, \, 0} \\ {0\, \, \, \, \, \, 1\, \, \, \, \, 0\, \, \, \, \, \, 0} \\ {0\, \, \, \, \, a\, \, \, \, \,\, b\,\, \, \, \, 0} \\ {c\, \, -c\, \, \, \, 0\, \, \, \, d} \end{array}\right) \, \, \, $ & $A=N \dotplus S, $ \newline $N=<e_3, e_4>,$ \newline $S=<e_1, e_2>$ & 2 & 0 & 1 \\ \hline
$As_{4}^{22} :$ & $ e_2e_2=e_2,$  $e_2e_4=e_4,$  \newline $e_3e_3=e_3,$   $e_3e_1=e_1,$ \newline $e_1e_2=e_1$ & $\left(\begin{array}{l} {a\, \, \, \, \, \, b\, \, -b\, \, \, \, 0} \\ {0\, \, \, \, \, \, 1\, \, \, \, \, \, \, 0\, \, \, \, \, 0} \\ {0\, \, \, \, \, 0\, \, \, \, \,\, \, \, 1\, \, \, \, \, 0} \\ {0\, \, \, \, \, \, c\, \, \, \, \, \, 0\, \, \, \, \, d} \end{array}\right)\, \, \, $ & $A=N \dotplus S, $ \newline $N=<e_1, e_4>,$ \newline $S=<e_2, e_3>$ & 2 & 1 & 0 \\ \hline
$As_{4}^{23}:$ & $ e_1e_1=e_4,$ $e_1e_4=-e_3,$ \newline $e_2e_1=e_3,$  $e_2e_2=e_3,$ \newline $e_4e_1=-e_3$ & $ \left(\begin{array}{l} {1\, \, \, \, \, \, 0\, \, \, \, \, 0\, \, \, \, \, \, \, \, \, \, \, \, \, 0} \\ {a\, \, \, \, \, \, 1\, \, \, \, \, 0\, \, \, \, \, \, \, \, \, \, \, \, \, \, 0} \\ {b\, \, \, \, \, c\, \, \, \, \, 1\, \, \, \, a(a+1)-2d} \\ {d\, \, \, \, \, a\, \, \, \, \, o\, \, \, \, \, \, \, \, \, \, \, \, \, \, 1}\end{array}\right)\, \, \, $& \textit{nilpotent} & 3 & 1 & 1 \\ \hline
$As_{4}^{24}:$ & $e_2e_2=e_2,$  $e_2e_3=e_3,$  \newline $e_2e_1=e_1,$  $e_4e_2=e_4,$ \newline $e_1e_2=e_1$ & $ \left(\begin{array}{l} {a\, \, \, \, 0\, \, \, \, 0\, \, \, \, 0} \\ {0\, \, \, \, 1\, \, \, \, 0\, \, \, \, 0} \\ {0\, \, \, \, b\, \, \, \, c\, \, \, \, 0} \\ {0\, \, \, \, d\, \, \, \, o\, \, \, \, e} \end{array}\right)\, \, \, $& $A=N \dotplus S, $ \newline $N=<e_1, e_3, e_4>,$ \newline $S=<e_2>$ & 3 & 1 & 1 \\ \hline
$As_{4}^{25}:$ & $e_1e_2=e_4,$  $e_1e_3=e_4,$ \newline $e_2e_1=-e_4,$   $e_2e_2=e_4,$ \newline $e_3e_1=e_4$ & $ \left(\begin{array}{l} {\, \, \, a\, \, \, \, \, \, \, \, \, \, 0\, \, \, \, \, \, 0\, \, \, \, \, 0} \\ {\, \, \, b\, \, \, \, \, \, \, \, \, \, a\, \, \, \, \, \, 0\, \, \, \, \, 0} \\ {-\frac{b^{2}}{2a} \, \, -b\, \, \, \, a\, \, \, \, 0} \\ {\, \, \, \, c\, \, \, \, \, \, \, \, \, d\, \, \, \, \, \, e\, \, \, \, a^{2} } \end{array}\right)\, \, \, $ & \textit{nilpotent} & 3 & 1 & 1 \\\hline
$As_{4}^{26} :$ & $e_1e_1=e_1,$  $e_1e_2=e_2,$ \newline $e_2e_1=e_2,$  $e_4e_1=e_4,$ \newline $e_3e_1=e_3$ & $\left(\begin{array}{l} {1\, \, \, \, \, 0\, \, \, \, \, 0\, \, \, \, \, 0} \\ {0\, \, \, \, a\, \, \, \, \, 0\, \, \, \, \, 0} \\ {b\, \, \, \, 0\, \, \, \, \, c\, \, \, \, \, d} \\ {e\, \, \, \, \, 0\, \, \, \, \, f\, \, \, g} \end{array}\right)\, \, \, $ & $A=N \dotplus S, $ \newline $N=<e_2, e_3, e_4>,$  \newline $S=<e_1>$ & 3 & 0 & 2 \\\hline

\end{tabular}

\begin{tabular}{|c|p{1.1in}|p{1.9in}|p{1.1in}|p{0.25in}|p{0.25in}|p{0.25in}|} \hline
$As_{4}^{27}:$ & $e_1e_1=e_1,$  $e_1e_2=e_2,$ \newline $e_1e_4=e_4,$  $e_1e_3=e_3,$ \newline $e_2e_1=e_2$ & $ \left(\begin{array}{l} {1\, \, \, \, \, 0\, \, \, \, \, 0\, \, \, \, \, 0} \\ {0\, \, \, \, a\, \, \, \, \, 0\, \, \, \, \, 0} \\ {b\, \, \, \, 0\, \, \, \, \, c\, \, \, \, \, d} \\ {e\, \, \, \, \, 0\, \, \, \, \, f\, \, \, g} \end{array}\right)\, \, \, $& $A=N \dotplus S, $ \newline $N=<e_2, e_3, e_4>,$  \newline $S=<e_1>$ & 3 & 2 & 0 \\\hline
$As_{4}^{28}(\alpha):$ & $e_1e_1=e_4,$  $e_1e_2=\alpha e_4,$ \newline $e_2e_1=-\alpha e_4,$ \newline $e_2e_2=e_4,$ $e_3e_3=e_4,$ & $\left(\begin{array}{l} {a\, \, \, b\, \, \,\, c\, \, \, \, \, \, \, \, \, d} \\ {e\, \, \, f\, \, \, g\, \, \, \, \, \, \,\, \, h} \\ {i\, \, \, j\, \, \,\, k\, \, \, \, \, \, \,\, \, \, l} \\ {0\, \, \, 0\, \, \, 0\, \, af-be} \end{array}\right)$ & \textit{nilpotent} & 3 & 1 & 1 \\ \hline
$As_{4}^{29}:$ & $e_1e_1=e_1,$  $e_1e_3=e_3,$  \newline $e_2e_2=e_2,$  $e_2e_4=e_4,$ \newline $e_3e_1=e_3,$  $e_4e_1=e_4$ & $\left(\begin{array}{l} {1\, \, \, \, \, \, 0\, \, \, \, \, 0\, \, \, \, \, 0} \\ {0\, \, \, \, \, 1\, \, \, \, \, \, 0\, \, \, \, \, 0} \\ {0\, \, \, \, \, 0\, \, \, \, \, a\, \, \, \, \, 0} \\ {b\, \, -b\, \, \, 0\, \, \, \, c} \end{array}\right)\, \, \, $ & \textit{unitary}, \newline $A=N \dotplus S, $ \newline $N=<e_3, e_4>,$  \newline $S=<e_1, e_2> $ & 3 & 0 & 0 \\ \hline
$As_{4}^{30} :$ & $e_1e_1=e_1,$  $e_1e_3=e_3,$ \newline $e_1e_4=e_4,$  $e_2e_2=e_2,$ \newline $e_3e_1=e_3,$ $e_4e_2=e_4$ & $\left(\begin{array}{l} {1\, \, \, \, \, \, 0\, \, \, \, \, 0\, \, \, \, \, 0} \\ {0\, \, \, \, \, 1\, \, \, \, \, \, 0\, \, \, \, \, 0} \\ {0\, \, \, \, \, 0\, \, \, \, \, a\, \, \, \, \, 0} \\ {b\, \, -b\, \, \, 0\, \, \, \, c} \end{array}\right)\, \, \, $ & \textit{unitary}, \newline $A=N \dotplus S, $ \newline $N=<e_3, e_4>,$  \newline $S=<e_1, e_2> $ & 3 & 0 & 0 \\ \hline
$As_{4}^{31} :$ & $e_1e_1=e_1,$  $e_1e_4=e_4,$ \newline $e_2e_2=e_2,$  $e_2e_3=e_3,$ \newline $e_3e_1=e_3,$ $e_4e_2=e_4$ & $\left(\begin{array}{l} {0\, \, \, \, \, 1\, \, \, \, \, \, 0\, \, \, \, \, 0} \\ {1\, \, \, \, \, \, 0\, \, \, \, \, \, 0\, \, \, \, \, 0} \\ {a\, \, -a\, \, \, \, 0\, \, \, \, b} \\ {c\, \, -c\, \, \, \, d\, \, \, \, 0} \end{array}\right)$, $\left(\begin{array}{l} {1\, \, \, \, \, 0\, \, \, \, \, \, 0\, \, \, \, \, 0} \\ {0\, \, \, \, \, \, 1\, \, \, \, \, \, 0\, \, \, \, \, 0} \\ {a\, \, -a\, \, \, \, b\, \, \, \, 0} \\ {c\, \, -c\, \, \, \, 0\, \, \, \, d} \end{array}\right)$ & \textit{unitary}, \newline $A=N \dotplus S, $ \newline $N=<e_3, e_4>,$  \newline $S=<e_1, e_2> $ & 2 & 0 & 0 \\ \hline
$As_{4}^{32} :$ & $e_2e_1=e_3,$  $e_3e_4=e_3$, \newline $e_4e_2=e_2,$  $e_4e_3=e_3,$ \newline $e_4e_4=e_4,$  $e_1e_4=e_1$ & $ \left(\begin{array}{l} {\, \, a\, \, \, \, \, \, 0\, \, \, \, \, \, 0\, \, \, -\frac{d}{b} \, } \\ {\, \, 0\, \, \, \, \, \, b\, \, \, \, \, 0\, \, \, \, \, \, \, \, c} \\ {-ac\, \, d\, \, \, \, ab\, \, \, \, \frac{cd}{b} } \\ {\, \, 0\, \, \, \, \, \, 0\, \, \, \, \, \, 0\, \, \, \, \, \, \, 1} \end{array}\right)\, \, \, $ & $A=N \dotplus S, $ \newline $N=<e_1, e_2, e_3>,$  \newline $S=<e_4> $ & 2 & 0 & 0 \\ \hline
$As_{4}^{33}:$ & $e_1e_1=e_2,$  $e_1e_2=e_3,$ \newline $e_1e_3=e_4,$  $e_2e_1=e_3,$ \newline $e_2e_2=e_4,$  $e_3e_1=e_4 $ & $\left(\begin{array}{l} {a\, \, \, \, \, \, \, \, \, \, \, 0\, \, \, \, \, \, \, \, \, \, \, \, 0\, \, \, \, \, \, \, \, \, \, 0} \\ {b\, \, \, \, \, \, \, \, \, \, a^{2} \, \, \, \, \, \, \, \, \, \, \, 0\, \, \, \, \, \, \, \, \, 0} \\ {c\, \, \, \, \, \, \, \, \, 2ab\, \, \, \, \, \, \, \, \, a^{3} \, \, \, \, \, \, \, \, 0} \\ {d\, \, \, 2ac+b^{2} \, \, 3a^{2} b\, \, \, \, a^{4} } \end{array}\right)\, \, $ & \textit{commutative,} \newline \textit{nilpotent} & 4 & 1 & 1 \\ \hline
$As_{4}^{34} :$ & $e_2e_2=e_2,$  $e_2e_3=e_3,$ \newline $e_2e_4=e_4,$  $e_2e_1=e_1,$ \newline $e_4e_2=e_4,$  $e_4e_3=e_1$ & $\left(\begin{array}{l} {ab\, \, \, \, c\, \, \, \, d\, \, \, \, ae} \\ {0\, \, \, \, \, 1\, \, \, \, \, 0\, \, \, \, \, 0} \\ {0\, \, \, \, \, e\, \, \, \, \, b\, \, \, \, \, 0} \\ {0\, \, \, \, \, 0\, \, \, \, 0\, \, \, \, \, a} \end{array}\right)\, \, \, $ & $A=N \dotplus S, $ \newline $N=<e_1, e_3, e_4>,$  \newline $S=<e_2>$  & 2 & 2 & 0 \\ \hline
$As_{4}^{35}:$ & $e_1e_2=e_1,$  $e_3e_2=e_3,$ \newline $e_2e_2=e_2,$  $e_2e_4=e_4,$ \newline $e_3e_4=e_1,$  $e_4e_2=e_4$ & $\left(\begin{array}{l} {ab\, \, \, \, c\, \, \, \, d\, \, \, \, ae} \\ {0\, \, \, \, \, 1\, \, \, \, \, 0\, \, \, \, \, 0} \\ {0\, \, \, \, \, e\, \, \, \, \, b\, \, \, \, \, 0} \\ {0\, \, \, \, \, 0\, \, \, \, 0\, \, \, \, \, a} \end{array}\right)\, \, \, $ &  $A=N \dotplus S, $ \newline $N=<e_1, e_3, e_4>,$  \newline $S=<e_2>$& 2 & 0 & 2 \\ \hline
$As_{4}^{36}:$ & $e_1e_1=e_1,$  $e_2e_2=e_2,$ \newline $e_2e_3=e_3,$  $e_2e_4=e_4,$ \newline $e_3e_1=e_3,$  $e_4e_1=e_4$ & $\left(\begin{array}{l} {1\, \, \, \, \, 0\, \, \, \, \, 0\, \, \, \, \, 0} \\ {0\, \, \, \, \, 1\, \, \, \, \, 0\, \, \, \, \, 0} \\ {a\, \, -a\, \, \, \, b\, \, \, \, \, c} \\ {d\, \, -d\, \, \, \, e\, \, \, f} \end{array}\right) \, \, \, $\textit{} & \textit{unitary}, \newline $A=N \dotplus S, $ \newline $N=<e_3, e_4>,$  \newline $S=<e_1, e_2> $ & 2 & 0 & 0 \\ \hline
$As_{4}^{37}:$ & $e_1e_1=e_1,$ $e_1e_2=e_2,$ \newline $e_1e_3=e_3,$  $e_2e_1=e_2,$ \newline $e_3e_1=e_3,$ $e_4e_1=e_4$ & $\left(\begin{array}{l} {1\, \, \, \, 0\, \, \, \, \, 0\, \, \, \, \, 0} \\ {0\, \, \, \, a\, \, \, \, b\, \, \, \, \, 0} \\ {0\, \, \, \, \, c\, \, \, \, d\, \, \, \, 0} \\ {e\, \, \, \, \, 0\, \, \, \, 0\, \, \, \, f} \end{array}\right)\, \, \, $ & $A=N \dotplus S, $ \newline $N=<e_2, e_3, e_4>,$  \newline $S=<e_1>$ & 3 & 0 & 1 \\ \hline
$As_{4}^{38}:$& $e_1e_1=e_1,$  $e_1e_2=e_2,$ \newline $e_1e_3=e_3,$  $e_1e_4=e_4,$ \newline $e_2e_1=e_2,$  $e_3e_1=e_3$ & $\left(\begin{array}{l} {1\, \, \, \, 0\, \, \, \, \, 0\, \, \, \, \, 0} \\ {0\, \, \, \, a\, \, \, \, b\, \, \, \, \, 0} \\ {0\, \, \, \, \, c\, \, \, \, d\, \, \, \, 0} \\ {e\, \, \, \, \, 0\, \, \, \, 0\, \, \, \, f} \end{array}\right)\, \, \, $ & $A=N \dotplus S, $ \newline $N=<e_2, e_3, e_4>,$  \newline $S=<e_1>$ & 3 & 1 & 0 \\ \hline
$As_{4}^{39}:$& $e_1e_1=e_1,$  $e_1e_2=e_2,$ \newline $e_1e_3=e_3,$  $e_2e_1=e_2,$ \newline $e_2e_2=e_3,$  $e_4e_1=e_4,$ \newline $ e_3e_1=e_3$ & $\left(\begin{array}{l} {1\, \, \, \, 0\, \, \, \, \, 0\, \, \, \, \, 0} \\ {0\, \, \, \, a\, \, \, \, 0\, \, \, \, \, 0} \\ {0\, \, \, \, b\, \, \, \, a^{2} \, \, 0} \\ {c\, \, \, \, \, 0\, \, \, \, 0\, \, \, \, d} \end{array}\right)\, \, \, $& $A=N \dotplus S, $ \newline $N=<e_2, e_3, e_4>,$  \newline $S=<e_1>$ & 3 & 0 & 1 \\ \hline
\end{tabular}

\begin{tabular}{|c|p{1.2in}|p{1.4in}|p{1.1in}|p{0.25in}|p{0.25in}|p{0.25in}|} \hline
$As_{4}^{40}:$ & $e_1e_1=e_1,$  $e_1e_2=e_2,$ \newline $e_1e_3=e_3,$   $e_1e_4=e_4,$ \newline $e_2e_1=e_2,$   $e_3e_1=e_3,$ \newline $e_4e_1=e_4$ & $\left(\begin{array}{l} {1\, \, \, \, 0\, \, \, \, \, 0\, \, \, \, \, 0} \\ {0\, \, \, \, a\, \, \, \, b\, \, \, \, \, c} \\ {0\, \, \, \, \, d\, \, \, \, e\, \, \, \, f} \\ {0\, \, \, \, \, i\, \, \, \, j\, \, \, \, k} \end{array}\right)\, \, \, $& \textit{commutative,} \newline \textit{unitary}, \newline $A=N \dotplus S, $ \newline $N=<e_2, e_3, e_4>,$  \newline $S=<e_1>$ & 4 & 0 & 0 \\ \hline
$As_{4}^{41}:$ & $e_1e_1=e_1,$  $e_1e_2=e_2,$ \newline $e_1e_4=e_4,$  $e_1e_3=e_3,$\newline $e_2e_1=e_2,$  $e_2e_2=e_3,$\newline $ e_3e_1=e_3$ & $ \left(\begin{array}{l} {1\, \, \, \, 0\, \, \, \, \, 0\, \, \, \, \, 0} \\ {0\, \, \, \, a\, \, \, \, 0\, \, \, \, \, 0} \\ {0\, \, \, \, b\, \, \, \, a^{2} \, \, 0} \\ {c\, \, \, \, \, 0\, \, \, \, 0\, \, \, \, d} \end{array}\right)\, \, $& $A=N \dotplus S, $ \newline $N=<e_2, e_3, e_4>,$  \newline $S=<e_1>$ & 3 & 1 & 0 \\ \hline
$As_{4}^{42}:$ & $e_1e_1=e_1,$  $e_1e_2=e_2,$ \newline $e_2e_3=e_1,$  $e_2e_4=e_2,$ \newline $e_3e_1=e_3,$  $e_3e_2=e_4,$ \newline $e_4e_3=e_3,$ $e_4e_4=e_4$ & $\left(\begin{array}{l} {\, \, \, \, 1\, \, \, \, \, \, \, \, \, \, \, a\, \, \, \, \, \, \, \, 0\, \, \, \, \, 0} \\ {\, \, \, \, 0\, \, \, \, \, \, \, \, \, \, \frac{\, 1}{b} \, \, \, \, \, \, \, \, 0\, \, \, \, \, 0} \\ {-ab\, \, -a^{2} b\, \, \, b\, \, \, \, ab} \\ {\, \, \, \, 0\, \, \, \, \, \, \, -a\, \, \, \, \, \, 0\, \, \, \, \, 1} \end{array}\right)\, \, \, $& \textit{unitary}, \newline $A=N \dotplus S, $ \newline $N=<e_2, e_3>,$  \newline $S=<e_1, e_4>$ & 2 & 0 & 0 \\ \hline
$As_{4}^{43}:$\textit{} & $e_1e_1=e_1,$  $e_1e_2=e_2,$ \newline $e_1e_3=e_3,$   $e_1e_4=e_4,$ \newline $e_2e_1=e_2,$ $e_2e_2=e_4,$ \newline $e_3e_1=e_3,$  $e_4e_1=e_4$ & $\left(\begin{array}{l} {1\, \, \, \, 0\, \, \, \, \, 0\, \, \, \, \, 0} \\ {0\, \, \, \, a\, \, \, \, 0\, \, \, \, \, 0} \\ {0\, \, \, \, b\, \, \, \, c\, \, \, \, \, 0} \\ {0\, \, \, \, d\, \, \, \, e\, \, \, \, a^{2} } \end{array}\right)\, \, \, $ & \textit{commutative,} \newline \textit {unitary}, \newline  $A=N \dotplus S, $ \newline $N=<e_2, e_3, e_4>,$  \newline $S=<e_1>$ & 4 & 0 & 0 \\ \hline
$As_{4}^{44}(\alpha ):$ & $e_1e_1=e_1,$  $e_1e_2=e_2,$ \newline $e_1e_3=e_3,$  $e_1e_4=e_4,$ \newline $e_2e_1=e_2,$ $e_2e_3=\alpha e_4$,\newline  $e_3e_1=e_3,$  $e_3e_2=e_4,$ \newline$e_4e_1=e_4$ &
$\left(\begin{array}{l} {1\, \, \, \, 0\, \, \, \, \, 0\, \, \, \, \, \, 0} \\ {0\, \, \, \, 0\, \, \, \, a\, \, \, \, \, \, \, 0} \\ {0\, \, \, \, b\, \, \, \, 0\, \, \, \, \, \, \,0} \\ {0\, \, \, \, c\, \, \, \, d\, \, \, \, \, ab} \end{array}\right)$,
$\left(\begin{array}{l} {1\, \, \, \, 0\, \, \, \, \, 0\, \, \, \, 0} \\ {0\, \, \, \, a\, \, \, \, 0\, \, \, \, \, 0} \\ {0\, \, \, \, 0\, \, \, \, b\, \, \, \, \, 0} \\ {0\, \, \, \, c\, \, \, \, d\, \, \, \, \, ab} \end{array}\right)$,
$\left(\begin{array}{l} {1\, \, \, 0\, \, \, 0\, \, \, \, \, \, \, 0} \\ {0\, \, \, a\, \, \, b\, \, \, \, \,\, \, 0} \\ {0\, \, \, c\, \, \, d\, \, \, \, \, \, \, 0} \\ {0\, \, \, e\, \, \, f\, \, ad-bc} \end{array}\right)$
 & \textit{unitary},\newline  $A=N \dotplus S, $ \newline $N=<e_2, e_3, e_4>,$  \newline $S=<e_1>$ & 3 & 0 & 0 \\ \hline
$As_{4}^{45} :$ & $e_1e_1=e_1,$  $e_1e_2=e_2,$ \newline $e_1e_3=e_3,$  $e_1e_4=e_4,$ \newline $e_2e_1=e_2,$ $e_2e_2=e_3,$ \newline $e_2e_3=e_4,$  $e_3e_1=e_3,$ \newline $e_3e_2=e_4,$ $e_4e_1=e_4$ & $\left(\begin{array}{l} {1\, \, \, \, 0\, \, \, \, \, \, \, \, 0\, \, \, \, \, \, \, \, \, 0} \\ {0\, \, \, \, a\, \, \, \, \, \, \, 0\, \, \, \, \, \, \, \, \, 0} \\ {0\, \, \, \, b\, \, \, \, \, \, \, a^{2} \, \, \, \, \, \, 0} \\ {0\, \, \, \, c\, \, \, \, 2a^{2} b\, \, \, a^{3} } \end{array}\right)\, \, \, $ & \textit{commutative,} \newline \textit{unitary}, \newline $A=N \dotplus S, $ \newline $N=<e_2, e_3, e_4>,$  \newline $S=<e_1>$ & 4 & 0 & 0 \\ \hline
$As_{4}^{46}:$ & $e_1e_1=e_1,$ $e_1e_2=e_2,$\newline $e_1e_3=e_3,$ $e_1e_4=e_4,$\newline $e_2e_1=e_2,$  $e_2e_2=-e_4,$ \newline $e_2e_3=-e_4,$ $e_3e_1=e_3,$ \newline $e_3e_2=e_4,$  $e_4e_1=e_4$ & $\left(\begin{array}{l} {1\, \, \, \, 0\, \, \, \, \, \, \, \, 0\, \, \, \, \, \, \, \, \, 0} \\ {0\, \, \, \, a\, \, \, \, \, \, \, 0\, \, \, \, \, \, \, \, \, 0} \\ {0\, \, \, \, b\, \, \, \, \, \, \, a^{2} \, \, \, \, \, \, 0} \\ {0\, \, \, \, c\, \, \, \, 2a^{2} b\, \, \, a^{3} } \end{array}\right)\, \, \, $& \textit{unitary}, \newline  $A=N \dotplus S, $ \newline $N=<e_2, e_3, e_4>,$  \newline $S=<e_1>$ & 3 & 0 & 0 \\ \hline
\end{tabular}\\

\end{footnotesize} 

where $a, b, c, d, e, f, g, h, i, j, k, l, \alpha \in \mathbb{C}.$

\end{document}